\documentclass[11pt]{article}
\usepackage{latexsym,amssymb,amsmath,amsfonts,amsthm}
\makeatletter \topmargin=0mm\headheight=0mm\headsep=0mm
\newcommand{\R}{{\mathbb R}}
\newcommand{\Z}{{\mathbb Z}}
\newcommand{\N}{{\mathbb N}}

\newcommand{\ben}{\begin{eqnarray}}
\newcommand{\een}{\end{eqnarray}}
\newcommand{\beno}{\begin{eqnarray*}}
\newcommand{\eeno}{\end{eqnarray*}}
\newcommand{\beq}{\begin{equation}}
\newcommand{\eeq}{\end{equation}}

\newcommand{\na}{\nabla}

\newcommand{\al}{\alpha}

\newtheorem{Theorem}{Theorem}[section]
\newtheorem{Lemma}[Theorem]{Lemma}
\newtheorem{Def}{Definition}[section]
\newtheorem{Remark}{Remark}[section]
\newtheorem{Proposition}{Proposition}[section]

\begin{document}
\title{\bf Existence theorem and blow-up criterion of the strong solutions to the two-fluid MHD equation  in $\R^3$}
\author {Qionglei Chen  and Changxing  Miao\\
{\small Institute of Applied Physics and Computational Mathematics,}\\
{\small P.O. Box 8009, Beijing 100088, P.R. China.}\\
{\small (chen\_{}qionglei@iapcm.ac.cn and   miao\_{}changxing@iapcm.ac.cn)}}
\date{}

\maketitle
 \vspace{-1.2in} \vspace{.9in} \vspace{0.2cm}

{\bf Abstract.} We first give the local well-posedness of strong
solutions to the Cauchy problem of the 3D two-fluid  MHD equations,
and then study the blow-up criterion of the strong solutions. By
means of the Fourier frequency localization and Bony's paraproduct
decomposition, it is proved that the strong solution $(u,b)$ can be
extended after $t=T$ if either $u\in L^q_T(\dot B^{0}_{p,\infty})$
with $\frac{2}{q}+\frac{3}{p}\le 1$ and $b\in L^1_T(\dot
B^{0}_{\infty,\infty})$  or $(\omega, J)\in L^q_T(\dot
B^{0}_{p,\infty})$ with $\frac{2}{q}+\frac{3}{p}\le 2$, where
$\omega(t)=\na\times  u $ denotes the  vorticity of the velocity
and $J=\na\times  b$  stands for  the current density.
 \vspace{0.2cm}

 \vspace{0.2cm}


\section{Introduction}
\setcounter{equation}{0}

We are concerned with the following   two-fluid magnetohydrodynamics
equations in $\R^3$: \ben\label{1.1} \left\{
\begin{aligned}\label{1.1}
&u_t-\nu\Delta u+u\cdot \nabla u-b\cdot \nabla b+\nabla
(p+\frac{1}{2}b^2)=0,\\
&b_t-\alpha\Delta b_t-\eta\Delta b+u\cdot \nabla b-b\cdot\nabla u+h\nabla\times(J\times b)=0 ,\\
&\nabla\cdot u=\nabla\cdot b=0,\\
&u(0,x)=u_0(x),\quad b(0,x)=b_0(x),
\end{aligned}\right.
\een where $x\in \R^3, t\ge 0$,  $\nu$, $\eta$, $\alpha$, $h$
stands for kinematic viscosity, the resistivity, the electron
inertia term and the Hall coefficient respectively,  $u$, $b$
describes the flow velocity vector and the magnetic field vector
respectively, and $J=\nabla \times b$ is the current density, $p$
is a scalar pressure, and  $u_0$ and $b_0$ are
 the given initial velocity and initial magnetic field
with $\nabla\cdot u_0=\nabla\cdot b_0=0$.
This model describes
some important physical phenomena.
In particular, for a plasma  composed of  two types
of fluids and formed by ions and electrons, this model can explain the phenomena of fast
magnetic reconnection  such as  in solar flares which cannot be characterized
appropriately by the one-fluid magnetohydrodynamics. It is
generally accepted now that  the two-fluid magnetohydrodynamics is more
complete than the classical  one-fluid magnetohydrodynamics(MHD)
model ( see  \cite{Bisk1,Bisk2,Priest} and references therein ).
This is the reason  why the two-fluid MHD equations are studied.

In general, the coefficient $\alpha$ is very small. Meanwhile, the Hall current term $h\nabla\times(J\times b)$ is also small
in dense plasmas,  so at large scales its effect is less important than that of the velocity.
Neglecting both of them, that is, formally letting $\alpha=h=0$,  the equations (1.1) reduce to the classical  MHD equations.
Further, if we also omit  the kinematic viscosity $\nu$ and  the resistivity $\eta$, that is, formally let $\alpha=h=\nu=\eta=0$,
we then obtain the classical ideal  MHD equations.
Both the MHD and  the ideal  MHD equations,  which are called one-fluid magnetohydrodynamics,
 have been studied extensively and are similar in many aspects to the Navier-Stokes equations and Euler equations, respectively.

It is well-known \cite{Ser} that  the classical  MHD equations are locally well-posed for any given initial
datum $u_0, b_0\in H^s(\mathbb R^3)$, $s\ge3$.
In the case of the two-fluid MHD equations,  N\'{u}\={n}ez \cite{Nun} has proved the  existence and uniqueness
of local solutions to the system for either Dirichlet or periodic boundary conditions. His result is

\vskip0.15cm

\noindent{\bf Theorem A.}\,   {\it If $u_0\in V$, $b_0\in D(A)$,  then there
exists an interval $[0,T]$ such that the two-fluid MHD equations
(\ref{1.1}) have a unique solution $(u,b)$ in $[0,T]$. Moreover
\begin{align}
&u\in C((0,T),V)\cap L^2((0,T),D(A)), \nonumber\\ & b\in
C((0,T),D(A)).\nonumber
\end{align}
Here in the Dirichlet case, $\Omega\subset\R^3$ is bounded and smooth, and
$$H=\{f\in L^2(\Omega)^3;\, \na\cdot f=0,\, f\cdot n|_{\partial\Omega}=0\}, \quad V=H^1_0(\Omega)^3\cap H, \quad
D(A)=H^2(\Omega)^3\cap V.$$
While, in the periodic case, $\Omega\subset\R^3$ is a box, and
\beno && H=\{f\in L^2(\Omega)^3;\, \int_{\Omega} f(x) dx=0,\, \na\cdot f=0,\, f\cdot n\,\, \mbox{antiperiodic  at  opposite  sides  of}\,\, \Omega\},\\
&& V=H^1(\Omega)^3\cap H, \quad
D(A)=H^2(\Omega)^3\cap V.\eeno}

The method of N\'{u}\={n}ez's proof seems not to apply to the Cauchy Problem \eqref{1.1} since
 the  Poincar\'{e} inequality plays a basic role in the  proof.
 The first purpose of this paper is to show the local well-posedness of strong solutions to the  equations  (\ref{1.1})
in $\R^3$ by Fourier localization together with Picard's method.

\begin{Theorem}\label{Th1}
If $(u_0, b_0)\in H^s(\R^3)\times H^{s+1}(\R^3)$ with $\nabla\cdot u_0=\nabla\cdot b_0=0$, $s\ge3$,  then there exists an
interval $[0,T]$ such that the two-fluid MHD equations (\ref{1.1}) have  a unique solution $(u,b)$ in
$C([0,T],(H^s\times H^{s+1})(\R^3))$. Moreover, $u$ satisfies
\beno
\begin{split}
u\in L^2([0,T],H^{s+1}(\R^3)).
\end{split}
\eeno
\end{Theorem}
Strong solutions we obtain here exist only locally.  In general, even for the classical MHD equations,
it is not known whether the smooth solution of the Cauchy problem
exists for all time though Duvaut and Lions\cite{Duv} constructed a class of global weak solutions.
An interesting question is  whether  smooth solutions will blow up at $t=T$
or,  in other words,  solutions can be extended to $[0,T')$  for $T'>T$  with the same regularity.
In particular, we want to obtain  conditions  under which the smooth solution  loses  its regularity at $t=T$
or the solution can be extended beyond $t=T$.

As we known,  for the 3D incompressible  Navier-Stokes equations, Giga\cite{Giga}
and Kozono-Taniuchi\cite{KozTa} obtained  criterions on extenstion
of strong solutions, that is,  strong solutions  can be continued
beyond $t=T$ provided  one of  following conditions holds:\vspace{.10cm}

$(1)\,\, u\in L^q(0,T; L^p(\R^3))$ for $\frac{2}{q}+\frac{3}{p}\le1, 3<p\le\infty,$\vspace{.1cm}

$(2)\,\, u\in L^2(\varepsilon_0, T; BMO(\R^3))$,\vspace{.1cm}

$(3)\,\, \nabla\times u\in L^1(\varepsilon_0, T; BMO(\R^3))$,\vspace{.15cm}\\
for $0\le \varepsilon_0<T$, where  $BMO$ is the space of  bounded mean oscillation functions.
On the other hand, many authors (see \cite{Bei,Koz} and references therein)
have studied  the regularity criterion for the weak solution such as:\vspace{.10cm}

$(4)\,\, \na u\in L^q(0,T; L^p(\R^3))$ for $\frac{2}{q}+\frac{3}{p}\le2$ and  $\frac32<p\le\infty$,\vspace{.1cm}

$(5)\,\, \na\times u\in L^q(0,T; \dot B^0_{p,\infty}(\R^3))$ for $\frac{2}{q}+\frac{3}{p}\le2$ and  $3\le p\le\infty$,\vspace{.1cm} \\
where $\dot B^0_{p,\infty}$ is homogeneous Besov space (see Section 2).

Caflisch, Klapper and Steele\cite{Caf}
extended the well-known  result  of  Beale-Kato-Majda\cite{Bea} for incompressible  Euler equations
to the cases of  the 3D ideal MHD equations.
Precisely, they showed that if the
smooth solution $(u,b)$ satisfies the condition
\ben\label{1.2}
\int_0^T\big(\|\nabla\times u(t)\|_{\infty}+\|\nabla\times b(t)\|_{\infty}\big)dt<\infty,
\een
then the solution $(u,b)$ can be extended beyond $t=T$. In other words,
let $[0,T)$ be the maximal existence time interval for  the smooth solution $(u, b)$
to the 3D ideal MHD equations.  Then  $(u, b)$ blows up at $T$ iff
\ben\label{1.3}
\lim_{\varepsilon\rightarrow0}\int_{\varepsilon_0}^{T-\varepsilon}
\big(\|\nabla\times u(t)\|_{\infty}+\|\nabla\times b(t)\|_{\infty}\big)dt=\infty, \qquad  \forall \  0\le\varepsilon_0<T.
\een
Recently,  the blow-up criterion \eqref{1.3} has been extended to
mixed time-space Besov spaces by the Fourier localization method (see \cite{CCM,ZL}).
For the classical MHD equations, Wu\cite{Wu0} showed that if the velocity and the magnetic field $(u,b)$ satisfy
\begin{align}\label{1.4}
\int_0^T(\|\na u(t)\|_2^4+\|\na b(t)\|_2^4)dt<\infty
\end{align}
or
\begin{align}\label{1.5}
\int_0^T(\|u(t)\|_\infty^2+\|b(t)\|_\infty^2dt)<\infty,
\end{align}
then the solution remains smooth. Later,
He and  Xin \cite{He} or  Zhou \cite{Zhou} obtained some integrability condition of
 the velocity $u$ alone, or  of the gradient  of the velocity   $\na u$ alone to characterize the
 regularity criterion for solutions to the classical MHD equations:
\begin{align}\label{1.6}
\int_0^T \|u(t)\|^q_{p}dt<\infty,\quad\,  \frac{2}{q}+\frac{3}{p}\le1\quad 3<p\le\infty
\end{align}
or
\begin{align}\label{1.7}
\int_0^T\|\na u(t)\|_{p}^qdt<\infty,\quad\,  \frac{2}{q}+\frac{3}{p}\le2 \quad \frac32<p\le\infty.
\end{align}
Other relevant results can be found in \cite{He,Wu,Zhou}.

As  mentioned above, there are similarities between the one-fluid MHD equations
and the Navier-Stokes equations. It is natural
to ask  whether similar  results hold for the two-fluid MHD equations.
The second purpose  of this paper is  to derive a similar blow-up criterion for the strong solution to
the  3D two-fluid MHD equations.
However, it seems to be difficult to obtain blow-up criterions using only the velocity $u$ like
  (\ref{1.6}) and  (\ref{1.7}).
 Roughly speaking, for  the classical MHD equations (i.e. $\alpha=h=0$), for given $u$, the magnetic induction equation is linear,
 so $b$ can be dominated by $\na u$ in some ways.
However,  for the 3D two-fluid MHD equations (i.e. $\alpha, h\neq0$),  the magnetic
induction equation is nonlinear with the nonlinear current
term $\na\times(J\times b)$,
so  the ``good'' term $-\alpha\Delta b_t-\eta\Delta b$ cannot compensate for the ``bad" effect caused by this nonlinear term.
This is why our blow-up criterion is given in terms of  both the velocity and the magnetic field.
We expect to establish a blow-up condition either on vorticity of $(u,b)$ or on $(u,b)$ in terms of Besov spaces  as in
\cite{Koz},
whose proof is based on the logarithmic Sobolev inequalities. However, in order to obtain the blow-up
criterion on $(u,b)$ itself, it seems that the logarithmic Sobolev inequalities do not work.
More precisely,  from the logarithmic Sobolev inequalities, one can deduce the following estimate of the solutions

$$f(t)\le C\exp\bigg(\int_0^tg(t')(\log f(t'))^kdt'\bigg) $$
for some $k>1$, which does not imply that  $f(t)$ will blow up in the finite time.
 To overcome this difficulty, we
make use of  the method of Fourier frequency localization and
Bony's paraproduct decomposition which enable us to obtain more precise nonlinear estimates. On the other hand,
for the blow-up condition on  the vorticity of $(u,b)$,
 our method gives a priori estimate  with one exponential growth, but the  logarithmic Sobolev inequalities
 only give  a priori estimate with a double exponential
growth. We now state our blow-up result.
\begin{Theorem}\label{Th2}
Assume that the initial solenoidal velocity and magnetic field $u_0\in H^s(\R^3)$, $b_0\in H^{s+1}(\R^3)$, $s\ge 3$.
Suppose that $(u,b)\in C([0,T],(H^s\times H^{s+1})(\R^3))$ is the strong solution to \eqref{1.1}.
If either
\begin{equation} \label{1.8}
\int_0^T(\|u(t)\|_{\dot B^{0}_{p,\infty}}^q+\|b(t)\|_{\dot B^{0}_{\infty,\infty}})dt<\infty \quad \mbox{with}\quad
 \frac{2}{q}+\frac{3}{p}\le 1, \quad 3< p\le\infty,
\end{equation}
or
\begin{equation} \label{1.9}
\int_0^T(\|\omega(t)\|_{\dot B^{0}_{p,\infty}}^q+\|J(t)\|_{\dot B^{0}_{p,\infty}}^q)dt<\infty \quad \mbox{with}\quad
 \frac{2}{q}+\frac{3}{p}\le 2, \quad 3\le  p\le\infty,
\end{equation}
then the solution $(u,b)$ can be extended beyond $t=T$. In other words, the solution blows up at $t=T$ iff
either
\beq\label{1.10}
\lim_{\varepsilon\rightarrow0}\int_{\varepsilon_0}^{T-\varepsilon}\big(\|u(t)\|^q_{\dot B^{0}_{p,\infty}}
+\|b(t)\|_{\dot B^{0}_{\infty,\infty}}\big)dt=\infty
\quad \mbox{with}\quad \frac{2}{q}+\frac{3}{p}\le 1,\quad 3<p\le\infty,
\eeq
or
\beq\label{1.11}
\lim_{\varepsilon\rightarrow0}\int_{\varepsilon_0}^{T-\varepsilon}
(\|\omega(t)\|_{\dot B^{0}_{p,\infty}}^q+\|J(t)\|_{\dot B^{0}_{p,\infty}}^q)dt=\infty \quad \mbox{with}\quad
 \frac{2}{q}+\frac{3}{p}\le 2, \quad 3\leq p\le\infty,
\eeq
where $\omega(t)=\na\times  u $ denotes the  vorticity of the velocity and $J=\na\times  b$ denotes the current density.
\end{Theorem}
\begin{Remark}
When $\alpha=h=0$, it is  known that if $(u,b)$ solves \eqref{1.1} then so does the pair of family $(u_\lambda, b_\lambda)$
for all $\lambda>0$, where $u_\lambda=\lambda u(\lambda x,\lambda^2t)$, $b_\lambda=\lambda b(\lambda x,\lambda^2t)$. Moreover,
$\|u_\lambda\|_{L^q(\R^+; L^p(\R^3))}=\|u\|_{L^q(\R^+; L^p(\R^3))}$ holds if and only if $\frac{2}{q}+\frac{3}{p}=1$.
However, in the case of
$\alpha, h\neq0$, the second equation of (\ref{1.1}) does not have such scaling invariance under the transformation
$(u,b)\mapsto(u_\lambda, b_\lambda)$. This is why we cannot set up a similar blow-up condition for the magnetic field $b$
as in the 3D MHD equations.
\end{Remark}
\begin{Remark}  In the conditions \eqref{1.9} and \eqref{1.11}, the integrability range of
 $\omega$ can be $\frac{3}{2}<p\le\infty$ by the Sobolev embedding theorem. On the other hand,
 by means of the H\"{o}lder inequality
$$\|b\|_{L^1_T({\dot B}^0_{\infty,\infty})}\le T^{1-\frac1{\tilde q}}\|b\|_{L^{\tilde q}_T({\dot B}^0_{\infty,\infty})},
\quad 1\le {\tilde q}\le\infty,$$
  the condition on $b$  in  \eqref{1.8} can  be extended to
 $ b\in L^{\tilde q}_T(\dot B^{0}_{\infty,\infty})$, $1\le {\tilde q}\le\infty$.  For the  case $p=\infty$,
 the two-fluid system seems
 to get  a benefit from  the term $-\alpha\Delta b_t$. For the classical MHD equations , we have the restriction
 condition  ${\tilde q}\ge 2$ in the case  $p=\infty$ ( see \cite{Wu0}).
\end{Remark}

\begin{Remark}
By means of the Sobolev embedding theorem $L^p\hookrightarrow {\dot B}^0_{p,\infty}$, $1\le p\le\infty$,
 the corresponding result to  Theorem \ref{Th2} can be obtained  in the framework  of Lebesgue spaces,
  that is, if either

\begin{equation*} \left\{ \begin{aligned}
        &u\in L^q_T( L^p) \quad \mbox{with}\quad
 \frac{2}{q}+\frac{3}{p}\le 1, \quad 3< p\le\infty, \\
                  & b\in L^1_T(L^{\infty}),
                          \end{aligned} \right.
                          \end{equation*}
or
\begin{equation*}  (\omega, J)\in L^q_T(L^p)\quad \mbox{with}\quad \frac{2}{q}+\frac{3}{p}\le
2,\quad 3\le p\le\infty, \end{equation*}
 then the solution $(u,b)$ of \eqref{1.1} can be
extended beyond $t=T$, where $L^q_T(X)$ denotes
$L^q((0,T);X)$.
\end{Remark}

\noindent{\bf Notation:} Throughout the paper, $C$ stands for   a generic  constant.
We  will use the notation $A\lesssim B$ to denote the relation  $A\le CB$ and
the notation $A\approx B$ to denote the relations  $A\lesssim B$ and $B\lesssim A$.
Further, $\|\cdot\|_{p}$ denotes the norm of the Lebesgue space $L^p$ and
$\|(f_1, f_2, \cdots, f_i)\|_{X}^a$ denotes $\|f_1\|_{X}^a+  \cdots+\|f_i\|_{X}^a$.

\section{Preliminaries}
\setcounter{equation}{0}

 Let us recall the Littlewood-Paley
decomposition.
Let ${\cal S}(\R^3)$ be
the Schwartz class of rapidly decreasing functions. Given $f\in
{\cal S}(\R^3)$, its Fourier transform ${\cal F}f=\hat f$ is
defined by
$$
\hat f(\xi)=(2\pi)^{-\frac{3}2}\int_{\R^3}e^{-ix\cdot \xi}f(x)dx.
$$
Choose two nonnegative radial functions $\chi$, $\varphi \in {\cal
S}(\R^3)$ supported respectively in ${\cal B}=\{\xi\in\R^3,\,
|\xi|\le\frac{4}{3}\}$ and ${\cal C}=\{\xi\in\R^3,\,
\frac{3}{4}\le|\xi|\le\frac{8}{3}\}$ such that
\beno
\chi(\xi)+\sum_{j\ge0}\varphi(2^{-j}\xi)=1,\quad\xi\in\R^3,\\
\sum_{j\in\Z}\varphi(2^{-j}\xi)=1,\quad\xi\in\R^3\backslash \{0\}.
\eeno
Set $\varphi_j(\xi)=\varphi(2^{-j}\xi)$ and
let $h={\cal F}^{-1}\varphi$ and $\tilde{h}={\cal F}^{-1}\chi$.
Define the frequency localization operators:
\beno
&&\Delta_jf=\varphi(2^{-j}D)f=2^{3j}\int_{\R^3}h(2^jy)f(x-y)dy, \\
&&S_jf=\sum_{k\le j-1}\Delta_kf=\chi(2^{-j}D)f=2^{3j}\int_{\R^3}\tilde{h}(2^jy)f(x-y)dy.
\eeno
Formally, $\Delta_j=S_{
j}-S_{j-1}$  is a frequency projection into the annulus
$\{|\xi|\approx 2^j\}$, and  $S_j$ is a frequency projection into the
ball $\{|\xi|\lesssim 2^j\}$. One easily verifies that
with the above choice of $\varphi$
\begin{eqnarray}\label{2.1}
\Delta_j\Delta_kf\equiv0\quad i\!f\quad|j-k|\ge 2\quad and
\quad \Delta_j(S_{k-1}f\Delta_k
f)\equiv0\quad i\!f\quad|j-k|\ge 5.
\end{eqnarray}
We now introduce the following definition of Besov spaces.
\begin{Def}\label{Def2.1}Let
$s\in \R, 1\le p,q\le\infty$. The homogenous Besov space $\dot
{B}^s_{p,q}$ is defined by
$$\dot {B}^s_{p,q}=\{f\in {\cal Z}'(\R^3):   \|f\|_{\dot
{B}^s_{p,q}}<\infty\},$$
\end{Def}
\noindent where
$$\|f\|_{\dot{B}^s_{p,q}}=\left\{\begin{array}{l}
\displaystyle\bigg(\sum_{j\in\Z}2^{jsq}\|\Delta_j f\|_p^q\bigg)^{\frac 1
q},\quad \hbox{for}\quad q<\infty,\\
\displaystyle\sup_{j\in \Z}2^{js}\|\Delta_jf\|_p, \quad \hbox{ for}
\quad q=\infty,
\end{array}\right.
$$
and ${\cal Z}'(\R^3)$ can be identified by the quotient
space ${\cal S}'/{\cal P}$ with the  space ${\cal
P}$ of polynomials.
\begin{Def}\label{Def2.1}Let
$s\in \R, 1\le p,q\le\infty$. The inhomogeneous Besov space $
{B}^s_{p,q}$ is defined by
$${B}^s_{p,q}=\{f\in {\cal S}'(\R^3): \|f\|_{
{B}^s_{p,q}}<\infty\},$$ \end{Def}
\noindent where
$$\|f\|_{{B}^s_{p,q}}=\left\{\begin{array}{l}
\displaystyle\bigg(\sum_{j\ge 0}2^{jsq}\|\Delta_j f\|_p^q\bigg)^{\frac 1
q}+\|S_0(f)\|_p,\quad \hbox{for}\quad q<\infty,\\
\displaystyle\sup_{j\ge 0}2^{js}\|\Delta_jf\|_p+\|S_0(f)\|_p, \quad
\hbox{ for} \quad q=\infty.
\end{array}\right.
$$
If $s>0$, then ${B}^s_{p,q}=L^p\cap\dot{B}^s_{p,q}$ and $\|f\|_{B^s_{p,q}}\approx\|f\|_{p}+\|f\|_{\dot{B}^s_{p,q}}.$
We refer to \cite{Ber,Tri} for  details.

Let us state some basic properties about the Besov spaces.
\begin{Proposition}\label{Prop2.1}
$(i)$ When $p=q=2$, the homogeneous Sobolev space $\dot H^s$ and $\dot B^2_{2,2}$
are equal and the two norms are equivalent:
$$\|f\|_{\dot H^s}\approx\|f\|_{\dot B^2_{2,2}}.$$
Similar properties hold for the the inhomogeneous Sobolev space $H^s$ and $B^2_{2,2}$.

$({ii})$ We have the equivalence of norms
$$\|D^k f\|_{\dot B^s_{p,q}}\approx\|f\|_{\dot B^{s+k}_{p,q}}, \quad \textrm{for}\quad  k\in \Z^+.$$
$({iii})$ Interpolation: for
$s_1, s_2\in\R$ and $\theta\in[0,1]$, one has $$\|f\|_{\dot B^{\theta s_1+(1-\theta)s_2}_{p,q}}\le
\|f\|^\theta_{\dot B^{s_1}_{p,q}}\|f\|^{(1-\theta)}_{\dot B^{s_2}_{p,q}}.$$
Similar interpolation inequality holds for inhomogeneous Besov spaces.\\
\end{Proposition}
The proofs of $({i})-({iii})$ are  standard and can be found in \cite{Ber,Tri}.\\

\section{Local existence and uniqueness}
\setcounter{equation}{0}

We now prove Theorem 1.1.

\noindent {\bf The proof of local existence}. It  involves the method of successive approximation.
Define the sequence $\{u^{(n)}, b^{(n)}\}_{n\in \N_0}$ by the following
linear system:
\ben\label{3.1} \left\{
\begin{aligned}\label{3.1}
&u^{(n+1)}_t-\nu\Delta u^{(n+1)}=-u^{(n)}\cdot \nabla u^{(n)}+b^{(n)}\cdot \nabla b^{(n)}-\nabla
(p^{(n)}+\frac{1}{2}{b^2}^{(n)}),\\
&b^{(n+1)}_t-\alpha\Delta b^{(n+1)}_t-\eta\Delta b^{(n+1)}=-u^{(n)}\cdot \nabla b^{(n)}+b^{(n)}\cdot\nabla u^{(n)}
-h\nabla\times(J^{(n)}\times b^{(n)}),\quad\\
&\nabla\cdot u^{(n+1)}=\nabla\cdot b^{(n+1)}=0,\\
&(u^{(n+1)}, b^{(n+1)})\big|_{t=0}=S_{n+2}(u_0,b_0),
\end{aligned}\right.
\een
where we set $(u^{(0)},b^{(0)})=(0,0)$, so $p^{(0)}=0$. We first derive the $L^2$ estimate of solutions.
By the divergence free condition, the embedding relation $H^s\hookrightarrow L^\infty$
and the $\epsilon$-Young inequality it easily to see that
\begin{align}\label{3.2}
&\frac{1}{2}\frac{d}{dt}\|\big(u^{(n+1)},\, b^{(n+1)},\, \alpha^\frac{1}{2}\na b^{(n+1)}\big)(t)\|_2^2+\nu\|\na
u^{(n+1)}(t)\|_2^2+\eta\|\na
b^{(n+1)}(t)\|_2^2\nonumber\\
\le &(\|u^{(n)}\|_2^2+\|b^{(n)}\|_2^2+h\|J^{(n)}\|_2\|b^{(n)}\|_2)(\|\na u^{(n+1)}\|_\infty+\|\na b^{(n+1)}\|_\infty)
\nonumber\\
\le& \frac{\nu}2\|\na u^{(n+1)}\|_{H^s}+
\frac{\eta}2\|\na b^{(n+1)}\|_{H^s}+
C\|(u^{(n)},\, b^{(n)},\, \alpha^\frac{1}{2}\na b^{(n)})\|_2^4.\end{align}
Now we derive  the $\dot H^s$ estimate. Apply the operator $\Delta_k$ to the equations (\ref{3.1}),
multiply the first one by $\Delta_ku^{(n+1)}$ and
the second one by $\Delta_kb^{(n+1)}$, integrate by parts to get, on noting that
${\rm div}u^{(n+1)}={\rm div}b^{(n+1)}=0$ , that
\begin{align}\label{3.3} &\frac{1}{2}
\frac{d}{dt}\|\big(\Delta_k u^{(n+1)},\,\Delta_kb^{(n+1)},\,\alpha^\frac{1}{2}\na\Delta_k b^{(n+1)}\big)(t)\|_2^2+\nu\|\na\Delta_k
u^{(n+1)}(t)\|_2^2+\eta\|\na\Delta_k
b^{(n+1)}(t)\|_2^2\nonumber\\
=&\big<\Delta_k(u^{(n)}\otimes u^{(n)}), \na\Delta_ku^{(n+1)}\big>
-\big<\Delta_k(b^{(n)}\otimes b^{(n)}),\, \na\Delta_ku^{(n+1)}\big>
\nonumber\\
&+\big<\Delta_k(u^{(n)}\otimes b^{(n)}-b^{(n)}\otimes u^{(n)}), \,\na\Delta_kb^{(n+1)}\big>
-h\big<\Delta_k(J^{(n)}\otimes b^{(n)}),\, \na\times\Delta_kb^{(n+1)}\big>,
\end{align}
where $\big<\cdot\,,\,\cdot\big>$ stands for the inner product. Multiplying $2^{2ks}$ on both sides of (\ref{3.3})
and  summing up over $k\in\Z$
yield that
\begin{equation}\label{3.4}
\frac{1}{2}\frac{d}{dt}\|\big(u^{(n+1)}, b^{(n+1)},\alpha^\frac{1}{2}\na b^{(n+1)}\big)(t)\|_{\dot H^s}
+\nu\|\na u^{(n+1)}(t)\|_{\dot H^s}+\eta\|\na b^{(n+1)}(t)\|_{\dot H^s}
\le\sum_{i=1}^4\Pi_i,
\end{equation}
where
\begin{align}\sum_{i=1}^4\Pi_i
\triangleq &
\sum_{k\in\Z}2^{2ks}\|\Delta_k(u^{(n)}\otimes u^{(n)})\|_{2}\|\Delta_k\nabla u^{(n+1)}\|_{2}
+\sum_{k\in\Z}2^{2ks}\|\Delta_k(b^{(n)}\otimes b^{n})\|_{2}\|\Delta_k\nabla u^{(n+1)}\|_{2}\nonumber\\&
+\sum_{k\in\Z}2^{2ks}\|\Delta_k(u^{(n)}\otimes b^{(n)}+b^{(n)}\otimes u^{n})\|_{2}\|\Delta_k\nabla b^{(n+1)}\|_{2}
\nonumber\\
&+h\sum_{k\in\Z}2^{2ks}\|\Delta_k(J^{(n)}\otimes b^{(n)})\|_{2}\|\Delta_k\nabla b^{(n+1)}\|_{2}.
\nonumber\end{align}
Using the Schwartz inequality, Lemma \ref{LemmaA.2} and
the embedding results that $H^s\hookrightarrow \dot H^s$ and
 $H^s\hookrightarrow L^\infty$, we obtain that
\begin{align}\label{3.5}
\Pi_1(t)&\le \|u^{(n)}u^{(n)}\|_{\dot H^s}\|\nabla u^{(n+1)}\|_{\dot H^s}
\le C\|u^{(n)}\|_{L^\infty}\|u^{(n)}\|_{\dot H^s}\|\nabla u^{(n+1)}\|_{\dot H^s}\nonumber\\
&\le\|u^{(n)}\|^2_{H^s}\|\nabla u^{(n+1)}\|_{H^s}
\le \frac{\nu}{4}\|\nabla u^{(n+1)}\|_{H^s}^2+C_\nu\|u^{(n)}\|_{H^s}^4.
\end{align}
Similarly, we have
\begin{align}\label{3.6}
\sum_{i=2}^4\Pi_i(t)\le \frac{\nu}{4}\|\nabla u^{(n+1)}\|_{H^s}^2+\frac{\eta}{4}\|\nabla b^{(n+1)}\|_{H^s}^2+
C_{\al,\nu,\eta,h}(\|(u^{(n)},b^{(n)},\al^{\frac12}\na b^{(n)})\|_{H^s}^4.
\end{align}
Set $E^{(n)}_s(t)\triangleq \|(u^{(n)},b^{(n)},\al^{\frac12}\na b^{(n)})\|_{H^s}^2,$ $n\in\N_0$.
Adding (\ref{3.2}) and  (\ref{3.4}),  and
using (\ref{3.5}) and  (\ref{3.6}) yield that
\begin{align}
&\frac{d}{dt}E^{(n+1)}_s(t)+\nu\|\na u^{(n+1)}(t)\|^2_{H^s}+\eta\|\na b^{(n+1)}(t)\|^2_{H^s}
\le \widetilde{C}\|(u^{(n)},b^{(n)},\,\al^{\frac12}\na b^{(n)})\|_{H^s}^4,\nonumber
\end{align}
where $\widetilde{C}=C_{\al,\nu,\eta,h}$. Integrating the above inequality with respect to $t$ gives
\begin{align}
&\sup_{t\in[0,T]}E^{(n+1)}_s(t)+\int_0^T\nu\|\na u^{(n+1)}(t)\|^2_{H^s}+\eta\|\na b^{(n+1)}(t)\|^2_{H^s}dt\nonumber\\
\le&\|S_{n+2}(u_0 , b_0 , \al^{\frac12} \na b_0)\|_{H^s}^2+\widetilde{C}
\int_0^T\|u^{(n)}(t)\|^4_{H^s}+\|b^{(n)}(t)\|^4_{H^s}+\al^2\|\na b^{(n)}(t)\|^4_{H^s}dt\nonumber\\
\le& C_0\|(u_0 , b_0 , \al^{\frac12} \na b_0)\|_{H^s}^2+\widetilde{C}T\big(\sup_{t\in[0,T]}E^{(n)}_s(t)\big)^2.\nonumber
\end{align}
Thus,  by the standard induction argument, it follows  that
\begin{align}\label{3.7}
&\|(u^{(n+1)},\,b^{(n+1)},\,\al^{\frac12}\na b^{(n+1)})(t)\|_{L^\infty_T(H^s)}+\nu^{\frac12}\|\na u^{(n+1)}(t)\|_{L^2_T(H^s)}
+\eta^{\frac12}\|\na b^{(n+1)}(t)\|_{L^2_T(H^s)}
\nonumber\\
\le& 2C_0\|(u_0 , b_0 , \al^{\frac12} \na b_0)\|_{H^s}
\end{align}
for all $n\in\N_0$,  and for  $T\in[0, T_0]$, where we set
$$T_0=\frac{1}{4C_0\widetilde{C}\|(u_0 , b_0 , \al^{\frac12} \na b_0)\|_{H^s}^2}.$$
Next we  show that there exists a positive time $T_1(\le T)$ independent of $n$
such that $\{u^{(n)}, b^{(n)}\}$ is a Cauchy sequence in the space
$${\cal X}^{s-1}_{T_1}\triangleq\big\{(f,g,\al^{\frac12}\na g)\in L^\infty_{T_1}(H^{s-1}),\,
(\nu^{\frac12}\na f, \eta^{\frac12}\na g)\in L^2_{T_1}(H^{s-1})\big\}.$$
Let $\delta u^{(n+1)}=u^{(n+1)}-u^{(n)}$, $\delta b^{(n+1)}=b^{(n+1)}-b^{(n)}$,
$\delta p^{(n+1)}=p^{(n+1)}-p^{(n)}$ and $\delta (b^2)^{(n)}={b^{2}}^{(n)}-{b^2}^{(n-1)}$  satisfy that
\ben \left\{
\begin{aligned}\label{3.8}
&\delta u^{(n+1)}_t-\nu\Delta\delta u^{(n+1)}=F_1+F_2+\cdots+F_5,\\
&\delta b^{(n+1)}_t-\alpha\Delta\delta  b^{(n+1)}_t-\eta\Delta\delta  b^{(n+1)}=G_1+G_2+\cdots+G_6
,\\&(\delta u^{(n+1)},  \delta b^{(n+1)})\big|_{t=0}=\Delta_{n+1}(u_0,b_0),
\end{aligned}\right.
\een
where
\begin{align}\sum_{j=1}^5F_i\triangleq&-\delta u^{(n)}\cdot \nabla u^{(n)}-u^{(n-1)}\cdot \nabla \delta u^{(n)}+\delta b^{(n)}\cdot \nabla b^{(n)}
+b^{(n-1)}\cdot \nabla \delta b^{(n)}-\nabla
(\delta p^{(n)}+\frac{1}{2}{\delta b^2}^{(n)}),\nonumber\\
\sum_{j=1}^6G_i\triangleq&-\delta u^{(n)}\cdot \nabla b^{(n)}-u^{(n-1)}\cdot \nabla \delta b^{(n)}+\delta b^{(n)}\cdot\nabla u^{(n)}
+b^{(n-1)}\cdot \nabla \delta u^{(n)}\nonumber\\&
-h\nabla\times(\delta J^{(n)}\times b^{(n)})-h\nabla\times(J^{(n-1)}\times \delta b^{(n)}).\nonumber\end{align}
Applying the divergence free condition to $F_2$, $F_3$ and $F_5$ yields that
\begin{align}\label{3.9}
\Big|\big<\sum_{j=1}^5F_i,\delta u^{(n+1)}\big>\Big|\le &\|\delta u^{(n)}\|_2\|\na u^{(n)}\|_\infty\|\delta u^{(n+1)}\|_2
+\|u^{(n-1)}\|_\infty\|\delta u^{(n)}\|_2
\|\na\delta u^{(n+1)}\|_2\nonumber\\
 &+\|\delta b^{(n)}\|_2\|\na b^{(n)}\|_\infty\|\delta u^{(n+1)}\|_2+\|b^{(n-1)}\|_\infty\|\delta b^{(n)}\|_2
\|\na\delta u^{(n+1)}\|_2^2\nonumber\\
\le &{\nu}\|\na\delta u^{(n+1)}\|^2_2+C\big(\|(u^{(n)},b^{(n)})\|_{H^s}+\|(u^{(n-1)},b^{(n-1)})\|^2_{H^s}\big)\nonumber\\
&\times\|(\delta u^{(n)},\,\delta b^{(n)},\,\delta u^{(n+1)})\|_2^2.
\end{align}
Similarly, we have
\begin{align}\label{3.10}
\Big|\big<\sum_{j=1}^6G_i,\delta u^{(n+1)}\big>\Big|&\le C\|(u^{(n)},\,b^{(n)},\,u^{(n-1)},\,b^{(n-1)})\|_{H^s}\nonumber \\ &\quad\times
\big(\|(\delta u^{(n)},\,\delta b^{(n)},\,\delta b^{(n+1)},\,\al^{\frac12}\na\delta b^{(n)},\,\al^{\frac12}\na\delta b^{(n+1)})\|_2^2\big).
\end{align}
Set $\delta E^{(n)}(t)\triangleq\|(\delta u^{(n)}, \delta b^{(n)}, \al^{\frac12}\na\delta b^{(n)})\|_2^2$.
By the $L^2$ energy estimate
combined with (\ref{3.9}) and  (\ref{3.10}) it is derived that
\begin{align}
\frac{d}{dt}\delta E^{(n+1)}(t)+\nu\|\na\delta u^{(n+1)}(t)\|_2^2+\eta\|\na\delta b^{(n+1)}(t)\|_2^2\le
C_1(\delta E^{(n)}(t)+\delta E^{(n+1)}(t)),\nonumber
\end{align}
where $C_1=C_{\al,\nu,\eta,h,\|(u_0, h_0, \sqrt{\al}\na h_0)\|_{H^s}^2}$. Integrating the above
inequality with respect to $t$ gives
\begin{align}\label{3.11}
&\sup_{t\in[0,T]}\delta E^{(n+1)}(t)+\int_0^T\nu\|\na\delta u^{(n+1)}(t)\|_2^2+\eta\|\na\delta b^{(n+1)}(t)\|_2^2 dt
\nonumber\\
\le & C_22^{-2(n+1)s}\|(u_0, h_0, \sqrt{\alpha}\nabla h_0)\|_{H^s}^2+C_1T\sup_{t\in[0,T]}(\delta E^{(n)}(t)+\delta E^{(n+1)}(t)),
\end{align}
where  use has been made of the fact that
$$\|\Delta_{n+1}(u_0, b_0, \al^{\frac12}\na b_0)\|_2
\le C_22^{-(n+1)s}\|(u_0, h_0, \sqrt{\alpha}\nabla h_0)\|_{H^s}.$$
 Thus,  if $C_1T\le\frac14$.
then we  have
\begin{align}
&\|(\delta u^{(n+1)}, \delta b^{(n+1)}, \al^{\frac12}\na\delta b^{(n+1)})\|_{L^\infty_T(L^2)}
+\nu^{\frac12}\|\na\delta u^{(n+1)}(t)\|_{L^2_T(L^2)}+\eta^{\frac12}\|\na\delta b^{(n+1)}(t)\|_{L^2_T(L^2)}\nonumber\\
\le &2C_22^{-(n+1)s},\qquad n\in\N_0.\nonumber
\end{align}
This,  together with the interpolation $\|f\|_{H^{s-1}}\le\|f\|_{2}^{\frac1s}\|f\|_{H^{s}}^{1-\frac{1}s}$ and (\ref{3.7}),
implies that
\begin{align}\label{3.12}
&\|(\delta u^{(n+1)}, \delta b^{(n+1)}, \al^{\frac12}\na\delta b^{(n+1)})\|_{L^\infty_T(H^{s-1})}
+\nu^{\frac12}\|\na\delta u^{(n+1)}(t)\|_{L^2_T(H^{s-1})}+\eta^{\frac12}\|\na\delta b^{(n+1)}(t)\|_{L^2_T(H^{s-1})}\nonumber\\
 \le& 2C_2C_02^{-(n+1)}
\|(u_0,b_0,\al^{\frac12}\na b_0)\|_{H^s}.
\end{align}
By a standard argument,  it can been shown that for $T_1\le\min\{T_0,\frac1{4C_1}\}$,
the sequence $\{u^{(n)}, b^{(n)}\}$ converges to $(u,b)$ in ${\cal X}_{T_1}^{s-1}$ which is
a equation to the equation  (\ref{1.1}).
Moreover, $(u,b)$ satisfies that
\begin{align}\label{3.13}
&\|(u, b, \al^{\frac12}\na b)(t)\|_{L^\infty_{T_1}(H^s)}+\nu^{\frac12}\|\na u(t)\|_{L^2_{T_1}(H^s)}+\eta^{\frac12}\|\na b(t)\|_{L^2_{T_1}(H^s)}
\nonumber\\&\le 2C_0\|(u_0 , b_0 , \al^{\frac12} \na b_0)\|_{H^s}.
\end{align}
{\bf The proof of the uniqueness}. Suppose  $(u', b')\in L^\infty_{T}(H^s)$ is another solution to (\ref{1.1}). Let $\delta u=u-u'$ and
$\delta b=b-b'$. Then $(\delta \theta,\delta u)$ satisfies the
following equations
\ben \left\{
\begin{aligned}
&\delta u_t-\nu\Delta\delta u=-\delta u\cdot\na u-u'\cdot\na\delta u+\delta b\cdot\na b+b'\cdot\na\delta b-\na(\delta p+\frac12\delta b^2),\\
&\delta b_t-\alpha\Delta\delta  b_t-\eta\Delta\delta  b=-\delta u\cdot\na b-u'\cdot\na\delta b+\delta b\cdot\na u+b'\cdot\na\delta u
-h\na\times(\delta J\times b)\nonumber\\&\hspace{3.9cm}-h(\na\times(J'\times\delta b)).
\end{aligned}\right.
\een
By the divergence free condition and integrating by part, we obtain that
\begin{align}
&\frac12\frac{d}{dt}\|(\delta u,\,\delta b,\,\al^{\frac12}\na\delta b)\|_2^2+\nu\|\na\delta u\|_2^2+\eta\|\na\delta b\|_2^2\nonumber\\
&=-\big<\delta u\cdot\na u, \delta u\big>+\big<\delta b\cdot\na b, \delta u\big>-\big<\delta u\cdot\na b, \delta b\big>
+\big<\delta b\cdot\na u, \delta b\big>-h\big<J'\times\delta b, \na\times\delta b\big>
\nonumber\\&\le\|\delta u\|_{2}^2\|\na u\|_\infty+2\|\delta b\|_2\|\delta u\|_2\|\na b\|_\infty+\|\delta b\|_{2}^2\|\na u\|_\infty
+h\|J'\|_\infty\|\delta b\|_2\|\na\times\delta b\|_2
\nonumber\\&\le C\|(u, b, b')\|_{H^s}(\|\delta u\|_2^2+\|\delta b\|_2^2+\al\|\na\delta b\|_2^2).\nonumber
\end{align}
Thus we have
\begin{align}
\|(\delta u,\,\delta b,\,\al^{\frac12}\na\delta b)\|_2\le C_3T\|(\delta u,\,\delta b,\,\al^{\frac12}\na\delta b)\|_2,\nonumber
\end{align}
where $C_3=C_{\|(u_0, b_0, \sqrt{\al}\na b_0)\|_{H^s}}$. This implies
that for sufficiently small $T\le T_2$, $\|(\delta u,\,\delta b,\,\al^{\frac12}\na\delta b)\|_2\equiv0$. Then by a
standard argument shows that the  uniqueness of local solutions in $L^\infty_{T}(H^s)$. This completes the proof
of Theorem 1.1.\endproof

\section{Blow-up criterion}
\setcounter{equation}{0}

In this section, we prove Theorem 1.2 which establishes the  blow-up criterion for the smooth solution to \eqref{1.1}.
The proof is broken down into two cases.

{\bf Case I.\quad The proof  of blow-up criterion under condition \eqref{1.8}.}\quad
We first derive a priori estimate of the smooth solution to
(\ref{1.1}).  Arguing similarly as in deriving (\ref{3.3}), we get
\begin{align}\label{4.1} &\frac{1}{2}
\frac{d}{dt}(\|\Delta_k u\|_2^2+\|\Delta_k
b\|_2^2+\alpha\|\na\Delta_k b\|_2^2)+\nu\|\na\Delta_k
u\|_2^2+\eta\|\na\Delta_k
b\|_2^2\nonumber\\
=&-\big<\Delta_k(u\cdot\na u), \Delta_k u\big>+
\big<\Delta_k(b\cdot\na b), \Delta_k u
\big>-\big<\Delta_k(u\cdot\na b), \Delta_k b\big>\nonumber\\&\quad+
\big<\Delta_k(b\cdot\na
u), \Delta_kb\big>-h\big<\Delta_k(\na\times(J\times b)), \Delta_k b
\big>.\end{align}
Noting that $\int_{\R^3}(b\times\Delta_k J)\Delta_k J
dx=0$, it follows that
\begin{align}
&-\big<\Delta_k(\na\times(J\times b)), \Delta_k b
\big>=\big<\Delta_k(b\times J), \Delta_k (\na\times b)
\big>=\big<(\Delta_k(b\times J)-b\times\Delta_k
J), \Delta_k J\big>,\nonumber
\end{align}
Substituting this into \eqref{4.1} and making use of the fact that  ${\rm div}u={\rm div}b=0$,  we obtain by
integrating by parts that
\begin{align}\label{4.2}&\frac{1}{2}
\frac{d}{dt}(\|\Delta_k u\|_2^2+\|\Delta_k
b\|_2^2+\alpha\|\na\Delta_k b\|_2^2)+\nu\|\na\Delta_k
u\|_2^2+\eta\|\na\Delta_k
b\|_2^2\nonumber\\
=&-\big<(\Delta_k(u\cdot\na
u)-u\cdot\na\Delta_k u), \Delta_k u \big>+\big<(\Delta_k(b\cdot\na
b)-b\cdot\na\Delta_k b), \Delta_k u\big>\nonumber\\
&+ \big<(\Delta_k(b\cdot\na
u)-b\cdot\na\Delta_k u), \Delta_k
b\big>-\big<(\Delta_k(u\cdot\na
b)-u\cdot\na\Delta_k b), \Delta_k
b\big>\nonumber\\
&+h\big<(\Delta_k(b\times
J)-b\times\Delta_k J), \Delta_k J\big>. \end{align}
Write the commutator $[f,\,\Delta_k]\cdot\na g$ for  $f\cdot\na\Delta_k g-\Delta_k(f\cdot\na g)$,
multiply  both sides of (\ref{4.2}) by $2^{2ks}$, and
 sum  the resulting equation over $k\in\Z$ to deduce that
\begin{align}\label{4.3}
&\frac{d}{dt}\|(u, b, \alpha^{\frac12}\na b)(t)\|_{\dot H^s}+2\nu\|\na u(t)\|_{\dot H^s}+\eta\|\na b(t)\|_{\dot H^s}
\nonumber\\
\le &\sum_{k\in\Z}2^{2ks}\|[u,\Delta_k]\cdot\na u\|_{2}\|\Delta_ku\|_{2}
+\sum_{k\in\Z}2^{2ks}\big(\|[b,\Delta_k]\cdot\na b\|_{2}\|\Delta_ku\|_{2}+\|[b,\Delta_k]\cdot\na u\|_{2}\|\Delta_kb\|_{2}\big)
\nonumber\\
&+\sum_{k\in\Z}2^{2ks}\|[u,\Delta_k]\cdot\na b\|_{2}\|\Delta_kb\|_{2}
+h\sum_{k\in\Z}2^{2ks}\|[b\times,\,\Delta_k]J\|_{2}\|\Delta_kJ\|_{2}
\nonumber\\
\triangleq& I+II+III+IV.
\end{align}
Making use of the   Schwartz inequality, and applying  Lemma
\ref{LemmaA.3} with $\sigma=s-1$, $\sigma_1=\sigma_2=-1$ and $p_1=p_2=p$ to the commutator, it follows that for $3<p\le\infty$,
\begin{align}\label{4.4}
|I|&\le C\big\|2^{k(s-1)}\|[u,\,\Delta_k]\cdot\na
u\|_{2}\big\|_{\ell^2(\Z)}\|u\|_{\dot B^{s+1}_{2,2}}
\nonumber\\&\le C\|u\|_{\dot B^0_{p,\infty}}\|u\|_{\dot
B^{s+\frac{3}{p}}_{2,2}}\|u\|_{\dot B^{s+1}_{2,2}}\le C\|u\|_{\dot
B^0_{p,\infty}}\|u\|_{\dot H^{s}}^{1-\frac{3}{p}}\|\na u\|_{\dot
H^s}^{1+\frac{3}{p}}\nonumber\\&\le {\nu}\|\na u\|_{\dot
H^s}^2+C\|u\|_{\dot B^0_{p,\infty}}^{\frac{2p}{p-3}}\|u\|_{\dot
H^s}^2,
\end{align}
where the use has been made of  the equivalent norms of
$\dot B^{\sigma}_{2,2}$ and $\dot H^\sigma$ for $\sigma\in\R$ and the interpolation theorem in deriving
the third inequality , and of the Young inequality to obtain the last inequality.
Similarly, for $3<p\le\infty$ we have the estimates
\begin{align}\label{4.5}
|II|&\le C\big\|2^{ks}\|[b,\,\Delta_k]\cdot\na b\|_{2}\big\|_{\ell^2(\Z)}\|u\|_{\dot B^{s}_{2,2}}+
\big\|2^{k(s-1)}\|[b,\,\Delta_k]\cdot\na u\|_{2}\big\|_{\ell^2(\Z)}\|b\|_{\dot B^{s+1}_{2,2}}\nonumber\\&
\le C\|b\|_{\dot B^0_{\infty,\infty}}\|b\|_{\dot B^{s+1}_{2,2}}\|u\|_{\dot B^{s}_{2,2}}
+C(\|b\|_{\dot B^0_{\infty,\infty}}\|u\|_{\dot B^{s}_{2,2}}+
\|u\|_{\dot B^0_{p,\infty}}\|b\|_{\dot B^{s+\frac{3}{p}}_{2,2}})\|b\|_{\dot B^{s+1}_{2,2}}\nonumber\\&\le
C(\|u\|_{\dot B^0_{p,\infty}}\|b\|_{\dot H^{s}}^{1-\frac{3}{p}}\|\na b\|_{\dot
H^s}^{1+\frac{3}{p}}+\|b\|_{\dot B^0_{\infty,\infty}}\|u\|_{\dot H^{s}}\|\na b\|_{\dot H^s})
\nonumber\\&\le \frac{\eta}{2}\|\na b\|_{\dot H^s}^2+C(\|u\|_{\dot B^0_{p,\infty}}
^{\frac{2p}{p-3}}+\|b\|_{\dot B^0_{\infty,\infty}})
(\|u\|_{\dot H^s}^2+\|b\|_{\dot H^s}^2+\alpha\|\na b\|_{\dot H^s}^2),\end{align}
and
\begin{align}\label{4.6}
|III|&\le C\big\|2^{k(s-1)}\|[u,\,\Delta_k]\cdot\na b\|_{2}\big\|_{\ell^2(\Z)}\|b\|_{\dot B^{s+1}_{2,2}}
\le C\big(\|u\|_{\dot B^0_{p,\infty}}\|b\|_{\dot H^{s+\frac{3}{p}}}
+\|b\|_{\dot B^0_{\infty,\infty}}\|u\|_{\dot H^{s}}\big)\|b\|_{\dot H^{s+1}}
\nonumber\\&\le \frac{\eta}{2}\|\na b\|_{\dot H^s}^2+C(\|u\|_{\dot B^0_{p,\infty}}
^{\frac{2p}{p-3}}+\|b\|_{\dot B^0_{\infty,\infty}})(\|u\|_{\dot H^s}^2+\|b\|_{\dot H^s}^2+\alpha\|\na b\|_{\dot H^s}^2)
\end{align}
and
\begin{align}\label{4.7}
|IV|&\le Ch\big\|2^{ks}\|[b\times,\,\Delta_k]J\|_{2}\big\|_{\ell^2(\Z)}\|J\|_{\dot B^{s}_{2,2}}
\nonumber\\&\le Ch\|b\|_{\dot B^0_{\infty,\infty}}
\|b\|_{\dot B^{s+1}_{2,2}}\|J\|_{\dot B^{s}_{2,2}}\le
Ch\|b\|_{\dot B^0_{\infty,\infty}}
\|\nabla b\|_{\dot H^{s}}^2.\end{align}
Integrating \eqref{4.3} with respect to $t$ and using
(\ref{4.4})-(\ref{4.7}), we deduce that for $3<p\le\infty$,
\begin{align}
&\sup_{t\in[0,T]}\|(u,\,b,\,\alpha^{\frac12}\na b)(t)\|_{\dot H^s}^2
+\int_0^T(\nu\|\na u(t)\|^2_{\dot H^s}+\eta\|\na b(t)\|_{\dot H^s}^2)dt
\nonumber\\
\le& \|(u_0,\, b_0,\, \alpha^{\frac12}\na b_0)\|_{\dot H^s}^2
+C\int_0^T(\|u(t)\|_{\dot B^0_{p,\infty}}^{\frac{2p}{p-3}}+\|b(t)\|_{\dot B^0_{\infty,\infty}})
\|(u,\,b,\,\alpha^{\frac12}\na b)(t)\|_{\dot H^s}^2dt.\nonumber
\end{align}
Note that $$0<\frac{2p}{p-3}\le q \quad\hbox{if}\quad \frac{2}{q}+\frac{3}{p}\le 1.$$
Then,  the Gronwall inequality yields
\begin{align}\label{4.8}
&\sup_{t\in[0,T]}(\|u(t)\|_{\dot H^s}^2+\|b(t)\|_{\dot H^s}^2+\alpha\|\na b(t)\|_{\dot H^s}^2)
+\int_0^T\Big(\nu\|\na u(t)\|^2_{\dot H^s}+\eta\|\na b(t)\|_{\dot H^s}^2\Big)dt
\nonumber\\&\le C(\|u_0\|_{\dot H^s}^2+\|b_0\|_{\dot H^s}^2+\alpha\|\na b_0\|_{\dot H^s}^2)\exp
\bigg(\|u(t)\|_{L^q_T(\dot B^0_{p,\infty})}^{\frac{2p}{p-3}}
T^{\frac{p}{p-3}(1-\frac{2}{q}-\frac{3}{p})}+\|b(t)\|_{L^1_T(\dot B^0_{\infty,\infty})}\bigg).
\end{align}
On the other hand, by the energy estimate we have
\begin{align}\label{4.9}
&\sup_{t\in[0,T]}(\|u(t)\|_{L^2}^2+\|b(t)\|_{2}^2+\alpha\|\na b(t)\|_{2}^2)
+\int_0^T\Big(\nu\|\na u(t)\|^2_{2}+\eta\|\na b(t)\|_{2}^2\Big)dt
\nonumber\\&\le \|u_0\|_{2}^2+\|b_0\|_{2}^2+\alpha\|\na b_0\|^2_{2}.
\end{align}
Combining  (\ref{4.8}) and (\ref{4.9}), and  by the standard argument of continuation of local solutions,
it is easy to show that if (\ref{1.8}) holds, then the solution remains smooth.

{\bf Case II.\quad The proof  of blow-up criterion under condition \eqref{1.9}.}\quad
Let us return to (\ref{4.3}). Applying  the Schwartz inequality and  Lemma
\ref{LemmaA.3} with $\sigma=s-\frac{3}{2p}$, $\sigma_1=\sigma_2=0$ and $p_1=p_2=p$ to the commutator,
it follows on using  the equivalent norms of
$\dot B^{\sigma}_{2,2}$ and $\dot H^\sigma$ for $\sigma\in\R$, the interpolation theorem and Young inequality,
that for $\frac 3 2<p<\infty$,
\begin{align}\label{4.10}
|I|&\le C\big\|2^{k(s-\frac{3}{2p})}\|[u,\,\Delta_k]\cdot\na u\|_{2}\big\|_{\ell^2(\Z)}\|u\|_{\dot B^{s+\frac{3}{2p}}_{2,2}}
\nonumber\\&\le C\|\na u\|_{\dot B^0_{p,\infty}}\|u\|^2_{\dot B^{s+\frac{3}{2p}}_{2,2}}\le C\|\na
u\|_{\dot B^0_{p,\infty}}\|u\|_{\dot H^{s}}^{\frac{2p-3}{p}}\|\na
u\|_{\dot H^s}^{\frac{3}{p}}\nonumber\\&\le \frac{\nu}{2}\|\na
u\|_{\dot H^s}^2+C\|\na u\|_{\dot B^0_{p,\infty}}^{\frac{2p}{2p-3}}\|u\|_{\dot H^s}^2.
\end{align}
Similar arguments as in deriving \eqref{4.10} can be used to  get that
\begin{align}\label{4.11}
|II|+|III|\le& C\big\|2^{k(s-\frac{3}{2p})}\|[b,\,\Delta_k]\cdot\na b\|_{2}\big\|_{\ell^2(\Z)}\|u\|_{\dot B^{s+\frac{3}{2p}}_{2,2}}+
\big\|2^{k(s-\frac{3}{2p})}\big(\|[b,\,\Delta_k]\cdot\na u\|_{2}\nonumber\\
&+\|[u,\,\Delta_k]\cdot\na b\|_{2}\big)\big\|_{\ell^2(\Z)}
\|b\|_{\dot B^{s+\frac{3}{2p}}_{2,2}}\nonumber\\
\le& C(\|\na b\|_{\dot B^0_{p,\infty}}\|u\|_{\dot B^{s+\frac{3}{2p}}_{2,2}}+
\|\na u\|_{\dot B^0_{p,\infty}}\|b\|_{\dot B^{s+\frac{3}{2p}}_{2,2}})\|b\|_{\dot B^{s+\frac{3}{2p}}_{2,2}}\nonumber\\
\le& C(\|\na u\|_{\dot B^0_{p,\infty}}+\|\na
b\|_{\dot B^0_{p,\infty}})(\|b\|_{\dot H^{s}}^{2-\frac{3}{p}}\|\na
b\|_{\dot H^s}^{\frac{3}{p}}+\|b\|_{\dot H^{s}}^{1-\frac{3}{2p}}\|\na
b\|_{\dot H^s}^{\frac{3}{2p}}\|u\|_{\dot H^{s}}^{1-\frac{3}{2p}}\|\na
u\|_{\dot H^s}^{\frac{3}{2p}})
\nonumber\\
\le &\bigg(\frac{\nu}{2}\|\na
u\|_{\dot H^s}^2+\frac{\eta}{2}\|\na b\|_{\dot H^s}^2\bigg)+C(\|\na u\|_{\dot B^0_{p,\infty}}
^{\frac{2p}{2p-3}}+\|\na
b\|^{\frac{2p}{2p-3}}_{\dot B^0_{p,\infty}})(\|u\|_{\dot H^s}^2+\|b\|_{\dot H^s}^2)
\end{align}
and
\begin{align}\label{4.12}
|IV|&\le Ch\big\|2^{ks}\|[b\times,\,\Delta_k]J\|_{2}\big\|_{\ell^2(\Z)}\|J\|_{\dot B^{s}_{2,2}}
\nonumber\\&\le Ch\|\na b\|_{\dot B^0_{p,\infty}}
\|b\|_{\dot B^{s+\frac{3}{p}}_{2,2}}\|J\|_{\dot B^{s}_{2,2}}\le
Ch\|\na b\|_{\dot B^0_{p,\infty}}\|b\|_{\dot H^{s}}^{1-\frac{3}{p}}
\|\nabla b\|_{\dot H^{s}}^{\frac{3}{p}}\|J\|_{\dot H^{s}}\nonumber\\&\le
\frac{\eta}{2}\|\na b\|_{\dot H^s}^2+C\|\na
b\|^{\frac{2p}{2p-3}}_{\dot B^0_{p,\infty}}(\|b\|_{\dot H^s}^2+\alpha\|\na b\|_{\dot H^s}^2).\end{align}
Integrating \eqref{4.3} with respect to $t$ and utilizing (\ref{4.10})-(\ref{4.12}) lead to the result  that for $3\le p<\infty$,
\begin{align}\label{4.13}
&\sup_{t\in[0,T]}\|(u,\,b,\,\alpha^{\frac12}\na b)(t)\|_{\dot H^s}^2+\int_0^T(\nu\|\na u(t)\|^2_{\dot H^s}+\eta\|\na b(t)\|_{\dot H^s}^2)dt
\nonumber\\
\le& \|(u_0,\,b_0,\,\alpha^{\frac12}\na b_0)\|_{\dot H^s}^2
+C\int_0^T(\|\na u(t)\|_{\dot B^0_{p,\infty}}^{\frac{2p}{2p-3}}+\|\na b(t)\|_{\dot B^0_{p,\infty}}^{\frac{2p}{2p-3}})
\|(u,\,b,\,\alpha^{\frac12}\na b)(t)\|_{\dot H^s}^2dt.
\end{align}
On the other hand, by the  Biot Savart law (\cite{Maj1}) we have
$$\nabla u=(-\Delta)^{-1}\nabla\nabla\times \omega,\quad
\nabla b=(-\Delta)^{-1}\nabla\nabla\times J,$$ where
$\omega=\na\times u$, $J=\na\times b$.   It follows from the  boundedness of
singular integral operators  on   homogeneous Besov spaces that
\begin{align}\label{4.14}\|(\na u, \na b)\|_{\dot B^\sigma_{p,r}}\le C\|(\omega, J)\|_{\dot B^{\sigma}_{p,r}} \quad
\textrm{for}\quad \sigma\in \R,\,\, (p,q)\in [1,\infty]\times[1,\infty].\end{align}
Inserting (\ref{4.14})  into (\ref{4.13}), we get
\begin{align}
&\sup_{t\in[0,T]}\|(u,\,b,\,\alpha^{\frac12}\na b)(t)\|_{\dot H^s}^2+\int_0^T(\nu\|\na u(t)\|^2_{\dot H^s}+\eta\|\na b(t)\|_{\dot H^s}^2)dt
\nonumber\\
\le& \|(u_0,\,b_0,\,\alpha^{\frac12}\na b_0)\|_{\dot H^s}^2+C\int_0^T(\|\omega(t)\|_{\dot B^0_{p,\infty}}^{\frac{2p}{2p-3}}+\|J(t)\|_{\dot B^0_{p,\infty}}^{\frac{2p}{2p-3}})
\|(u,\,b,\,\alpha^{\frac12}\na b)(t)\|_{\dot H^s}^2dt.\nonumber
\end{align}
Note that
$$0<\frac{2p}{2p-3}\le q \quad\hbox{if}\quad \frac{2}{q}+\frac{3}{p}\le 2.$$
Then,  the Gronwall inequality implies that for $3\le p<\infty$,
\begin{align}\label{4.15}
&\sup_{t\in[0,T]}\|(u,\,b,\,\alpha^{\frac12}\na b)(t)\|_{\dot H^s}^2\le
\|(u_0,\,b_0,\,\alpha^{\frac12}\na b_0)\|_{\dot H^s}^2\exp
\bigg(C\int_0^T\big(\|(\omega, J)(t)\|_{\dot{B}^0_{p,\infty}}^{\frac{2p}{2p-3}}
\big)dt\bigg)\nonumber\\
\le&\|(u_0,\,b_0,\,\alpha^{\frac12}\na b_0)\|_{\dot H^s}^2\exp
\bigg(\|(\omega, J)(t)\|_{L^q_T(\dot{B}^0_{p,\infty})}^{\frac{2p}{2p-3}}
T^{\frac{p}{2p-3}(2-\frac{2}{q}-\frac{3}{p})}\bigg).
\end{align}

For the case $p=\infty$, we apply Lemma \ref{LemmaA.3}   with $p_1=p_2=\infty$,
$\sigma_1=\sigma_2=0$ to the commutator to obtain that
\begin{align}
|I|+|II|+|III|+|IV|\le C(\|\na u\|_{L^\infty}+\|\na b\|_{L^\infty})(\|u\|_{\dot H^s}^2+\|b\|_{\dot H^s}^2+\|\na b\|_{\dot H^s}^2).\nonumber
\end{align}
Using the above estimate along with (\ref{4.3})  it follows from the Gronwall inequality that
\begin{align}\label{4.16}
&\sup_{t\in[0,T]}\|(u,\,b,\,\alpha^{\frac12}\na b)(t)\|_{\dot H^s}^2\le
\|(u_0,\,b_0,\,\alpha^{\frac12}\na b_0)\|_{\dot H^s}^2\exp
\bigg(C\int_0^T\|(\nabla u, \nabla b)(t)\|_{L^\infty}dt\bigg).
\end{align}
The logarithmic Sobolev inequality ( see (2.2) in \cite{Koz} )  and (\ref{4.14}) allow us to get that
\begin{align}\label{4.17}&\|(\na u,\na b)\|_{L^\infty}\le C\bigg(1+\|(\na u,\na
b)\|_{\dot{B}^0_{\infty,\infty}}
\log(\|(u,b)\|_{{\dot H^s}}+e)\bigg)\nonumber\\
\le& C\bigg(1+\|(\omega, J)\|_{\dot{B}^0_{\infty,\infty}}
\log(\|(u,b)\|_{{\dot H^s}}+e)\bigg).\end{align}
Plugging  (\ref{4.17}) into (\ref{4.16}), and setting  $Z(t)\triangleq\log(\|(u(t),b(t),\alpha^{\frac12}\na b(t))\|_{\dot H^s}+e)$,
we  deduce that
\begin{align}
\sup_{t\in[0,T]}Z(t)\le Z(0)+CT+C\int_0^T\|(\omega, J)(t)\|_{\dot{B}^0_{\infty,\infty}}Z(t)dt.\nonumber
\end{align}
Then the Gronwall inequality yields that
\begin{align}
&\sup_{t\in[0,T]}Z(t)\le (Z(0)+CT)\exp\bigg(C\|(\omega, J)(t)\|_{L^q_T(\dot{B}^0_{\infty,\infty})}
T^{(1-\frac{1}{q})}\bigg)\nonumber.
\end{align}
This implies that
\begin{align}\label{4.18}
&\sup_{t\in[0,T]}(\|u(t)\|_{\dot H^s}+\|b(t)\|_{\dot H^s}+\alpha^{\frac12}\|\na b(t)\|_{\dot H^s})\nonumber\\
\le&  (\|u_0\|_{\dot H^s}+\|b_0\|_{\dot H^s}+\alpha^{\frac12}\|\na b_0\|_{\dot H^s}+e)
^{A(T)}\exp(CTA(T)),
\end{align}
where $$A(T)\triangleq \exp\Big(C\|(\omega, J)(t)\|_{L^q_T(\dot{B}^0_{\infty,\infty})}
T^{(1-\frac{1}{q})}\Big).$$
Combining \eqref{4.18}  with (\ref{4.9}) and using the standard argument of continuation of local solutions,
we easily  prove that if (\ref{1.9}) holds, the solution remains smooth.
The proof of Theorem 1.2 is thus complete.
\endproof\vspace{.2cm}

\section{Appendix}
\setcounter{equation}{0}

Let us recall the paradifferential calculus which enables us to define a generalized
product between distributions, which is continuous in many functional spaces where the usual
product does not make sense (see \cite{Bon}). The paraproduct between $u$ and $v$ is defined
by$$T_uv\triangleq\sum_{j\in\Z}S_{j-1}u\Delta_jv.$$
We then have the following formal decomposition:
\beq\label{A.1}
uv=T_uv+T_vu+R(u,v)
\eeq
with $$R(u,v)=\sum_{j\in\Z}\Delta_ju\widetilde{\Delta}_jv\quad\mbox{and}\quad
\widetilde{\Delta}_j=\Delta_{j-1}+\Delta_j+\Delta_{j+1}.$$
The decomposition (\ref{A.1}) is called  Bony's paraproduct decomposition.

We first introduce the well-known  Bernstein inequality which
will be used repeatedly in the proof of the commutator estimate.

\begin{Lemma}\label{LemmaA.1}
Let $k$ be in $\N$. Let $(R_1, R_2)$ satisfy $0<R_1<R_2$. There exists a constant $C$ depending
only on  $R_1, R_2, d$ such that for all $1\le p\le q\le\infty$ and $f\in L^p(\R^d)$
we have
\begin{align}
&\sup_{|\gamma|=k}\|\partial^\gamma f\|_{q}\le C^{k+1}\lambda^{{k}+d(\frac{1}{p}-\frac{1}{q})}\|f\|_{p}\quad
\mbox{if}\quad{\rm supp}\hat f\subset {\cal B}(0, R_1\lambda),
\label{A.2}\\
&C^{-k-1}\lambda^{k}\|f\|_{p}\le C\sup_{|\gamma|=k}\|\partial^\gamma f\|_{p}\quad\mbox{if}\quad
{\rm supp}\hat f\subset {\cal C}(0, R_1\lambda, R_2\lambda).\label{A.3}
\end{align}
\end{Lemma}
\noindent The proof can be found in \cite{Ch1}.

\begin{Lemma}\label{LemmaA.2}
Let $s>0$, $f, g\in L^\infty\cap \dot H^s$.  Then $fg\in L^\infty\cap \dot H^s$ and
\begin{align}
\|fg\|_{\dot H^s}\le C(\|f\|_\infty\|g\|_{\dot H^s}+\|g\|_\infty\|f\|_{\dot H^s}).\nonumber
\end{align}
\end{Lemma}
\noindent For a proof see \cite{Ch1}.

\begin{Lemma}\label{LemmaA.3}Let $1\le p_1, p_2\le\infty$, $\sigma>0$, $\frac{d}{p_i}-\sigma_i>0$($i=1,2$) and
assume that $\sigma-\sigma_2+\frac{d}{p_2}>0$.
Then the following inequality holds:
\begin{align}\label{A.4}
\bigg(\sum_{j\in\Z}2^{2j\sigma}\|[f,\,\Delta_j]\na g\|^2_{2}\bigg)^{\frac12}\lesssim
\|\na f\|_{\dot B^{\sigma_1}_{p_1,\infty}}\|g\|_{\dot B^{\sigma-\sigma_1+\frac{d}{p_1}}_{2,2}}+
\|\na g\|_{\dot B^{\sigma_2}_{p_2,\infty}}\|f\|_{\dot B^{\sigma-\sigma_2+\frac{d}{p_2}}_{2,2}}.
\end{align}
If $\sigma_1=0,\, p_1=\infty$, $\|\na f\|_{\dot B^{\sigma_1}_{p_1,\infty}}$
has  to be replaced by $\|\na f\|_{L^\infty}$,  and  if $\sigma_2=0,\, p_2=\infty$,
then   $\|\na g\|_{\dot B^{\sigma_2}_{p_2,\infty}}$
has  to be replaced by $\|\na g\|_{L^\infty}$.
In \eqref{A.4}
 $$[f,\,\Delta_j]\na g=f\Delta_j(\na g)-\Delta_j(f\na g)$$
\end{Lemma}
\noindent{\it Proof:}\,\, The proof is standard. By Bony's decomposition, we have
\begin{align}\label{A.5}
[f,\,\Delta_j]\na g=[T_f,\,\Delta_j]\na g+T'_{\Delta_j\na g}f-\Delta_jT_{\na g}f-\Delta_jR(f,\na g),
\end{align}
where $T'_uv$ stands for $T_uv+R(u,v)$.

By (\ref{2.1}), we have
\begin{align}
[T_f,\,\Delta_j]\na g=
\sum_{j'\sim j}2^{jd}\int_{\R^3}h(2^j(x-y))(S_{j'-1}f(x)-S_{j'-1}f(y))\Delta_{j'}\na g(y)dy, \nonumber
\end{align}
where $j'\sim j$ means that  $|j'-j|\le4$.
Since $\frac{d}{p_1}-\sigma_1>0$, by Lemma \ref{LemmaA.1} and  the H\"{o}lder inequality we infer that
\begin{align}\label{A.6}
\|S_{j'-1}\na f\|_{\infty}\le C\sum_{j''\le j'-2}2^{j''\sigma_1}\|\Delta_{j''}\na f\|_{{p_1}}
2^{j''(\frac{d}{p_1}-\sigma_1)}
\le C\|\na f\|_{\dot B^{\sigma_1}_{p_1,\infty}}2^{j'(\frac{d}{p_1}-\sigma_1)}.
\end{align}
so
\begin{align}
\|[T_f,\,\Delta_j]\na g\|_{2}&\le C2^{-j}\sum_{j'\sim j}\|2^j|\cdot|2^{jd}h(2^j\cdot)\|_{1}
\|S_{j'-1}\na f\|_{{\infty}}\|\Delta_{j'}\na g\|_{2}\nonumber\\
& \le C\|\na f\|_{\dot
B^{\sigma_1}_{p_1,\infty}}2^{-j}\sum_{j'\sim j}2^{j'(\frac{d}{p_1}-\sigma_1)}
\|\Delta_{j'}\na g\|_{2}.\nonumber
\end{align}
This, together with the convolution inequality for series, gives
\begin{align}\label{A.7}\big\|2^{\sigma j}\|[T_f,\,\Delta_j]\na g\|_{2}\big\|_{\ell^2(\Z)}&\le C\|\na f\|_{\dot B^{\sigma_1}_{p_1,\infty}}
\bigg\|\sum_{j'\sim j}2^{(j'-j)(1-\sigma)}2^{j'(\sigma-\sigma_1+\frac{d}{p_1})}\|\Delta_{j'}g\|_{2}\bigg\|_{\ell^2(\Z)}\nonumber\\
&\le C\|\na f\|_{\dot B^{\sigma_1}_{p_1,\infty}}\|g\|_{\dot B^{\sigma-\sigma_1+\frac{d}{p_1}}_{2,2}}.
\end{align}
Using the definition of $T'_{\Delta_j\na g}f$ and (\ref{2.1}), we can rewrite
\begin{align}
T'_{\Delta_j\na g}f=\sum_{j'\gtrsim j}\Delta_{j'}f  S_{j'+2}\Delta_j\na g,\nonumber
\end{align}
where $j'\gtrsim j$ means that  $j'\ge j-2$.
By Lemma \ref{LemmaA.1} and the H\"{o}lder inequality, it follows that
\begin{align}
\|S_{j'+2}\Delta_j\na g\|_{\infty}\le C2^{j\sigma_2}\|\Delta_j\na g\|_{{p_2}}2^{j(\frac{d}{p_2}-\sigma_2)}
\le C\|\na g\|_{\dot B^{\sigma_2}_{p_2,\infty}}2^{j(\frac{d}{p_2}-\sigma_2)}.\nonumber
\end{align}
Thus,  for $\sigma-\sigma_2+\frac{d}{p_2}>0$, the convolution inequality for series yields that
\begin{align}\label{A.8}\big\|2^{\sigma j}\|T'_{\Delta_j\na g}f\|_{2}\big\|_{\ell^2(\Z)}&\le C\|\na g\|_{\dot B^{\sigma_1}_{p_2,\infty}}
\bigg\|\sum_{j'\gtrsim j}2^{-(j'-j)(\sigma-\sigma_2+\frac{d}{p_2})}2^{j'(\sigma-\sigma_2+\frac{d}{p_2})}
\|\Delta_{j'}f\|_{2}\bigg\|_{\ell^2(\Z)}\nonumber\\
&\le C\|\na g\|_{\dot B^{\sigma_2}_{p_2,\infty}}\|f\|_{\dot B^{\sigma-\sigma_2+\frac{d}{p_2}}_{2,2}}.
\end{align}
Similarly as in deriving (\ref{A.6}), we can show that  for $\frac{d}{p_2}-\sigma_2>0$,
we get $\|S_{j'-1}\na g\|_{\infty}\le C\|\na g\|_{\dot B^{\sigma_2}_{p_2,\infty}}2^{j'(\frac{d}{p_2}-\sigma_2)}.$
This, together with the convolution inequality for series, implies that
\begin{align}\label{A.9}
&\big\|2^{j\sigma}\|\Delta_j(T_{\na g}f)\|_{2}\big\|_{\ell^2(\Z)}\le
\bigg\|2^{j\sigma}\sum_{j'\sim j}\|\Delta_j(\Delta_{j'}f    S_{j'-1}\na g)\|_{2}\bigg\|_{\ell^2(\Z)}\nonumber\\
\le& C\|\na g\|_{\dot B^{\sigma_2}_{p_2,\infty}}\bigg\|\sum_{j'\sim j}2^{j'(\frac{d}{p_2}-\sigma_2+\sigma)}
\|\Delta_{j'}f\|_{2}2^{(j-j')\sigma}\bigg\|_{\ell^2(\Z)}
\le C\|\na g\|_{\dot B^{\sigma_2}_{p_2,\infty}}\|f\|_{\dot B^{\sigma-\sigma_2+\frac{d}{p_2}}_{2,2}}.
\end{align}
Finally, for $\sigma>0$ we have
\begin{align}\label{A.10}
&\big\|2^{j\sigma}\|\Delta_jR(f,\na g)\|_{2}\big\|_{\ell^2(\Z)}
\le\bigg\|\sum_{j'\gtrsim j}2^{j\sigma}\|\Delta_j(\Delta_{j'}f  \widetilde{\Delta}_{j'}\na g)\|_{2}\bigg\|_{\ell^2(\Z)}\nonumber\\
\le& C\|\na g\|_{\dot B^{\sigma_2}_{p_2,\infty}}\bigg\|\sum_{j'\gtrsim j}\|\Delta_{j'}f\|_{2}2^{j'(\frac{d}{p_2}-\sigma_2+\sigma)}
2^{(j-j')\sigma}\bigg\|_{\ell^2(\Z)}\le C\|\na g\|_{\dot B^{\sigma_2}_{p_2,\infty}}\|f\|_{\dot B^{\sigma-\sigma_2+\frac{d}{p_2}}_{2,2}}.
\end{align}

Combining (\ref{A.7})-(\ref{A.10}) gives  the desired inequality (\ref{A.4}).\endproof
\vspace{.2cm}
\textbf{Acknowledgements} The authors  would like  to thank  Professors G.
Ponce,   Z. Xin and  B. Zhang   so much for  their  helpful discussion and
suggestions. The authors are also deeply grateful to the referee
for his  valuable advices  which helped improve the paper greatly.  Q. Chen and C. Miao  were
supported in part by the NSF of China (No.10571016) and  C. Miao  was supported in part by
 The Institute of Mathematical Sciences, The Chinese
University of Hong Kong.


\begin{thebibliography}{50}



\bibitem{Bea} Beale J. T.,   Kato T.  and Majda A. J.,  {\it Remarks on the breakdown of
smooth solutions for the $3$-D Euler equations}, Comm. Math.
Phys., 94 (1984), 61-66.


\bibitem{Bei} Beirao da Veiga, H., {\it A new regularity class for the Navier-Stokes equations in $\R^n$},
Chin. Ann. Math., 16B (1995), 407-412.

\bibitem{Ber} Bergh J.   and   L\"{o}fstrom J., {\it Interpolation spaces, An
              Introduction},  New York: Springer-Verlag, 1976.

\bibitem{Bon}  Bony J.-M., {\it Calcul symbolique et propagation des singularit\'{e}s pour les \'{e}quations aux d\'{e}riv\'{e}es partielles
              non lin\'{e}aires}, Ann. Sci. \'{E}cole Norm. Sup., 14 (1981), 209-246.

\bibitem{Bisk1}  Biskamp D., {\it Nonlinear magnetohydrodynamics},  Cambridge University Press, Cambridge, UK, 1993.

\bibitem{Bisk2}  Biskamp D., {\it Magnetic reconnection in Plasmas },  Cambridge University Press, Cambridge, UK, 2000.



\bibitem{Caf}  Caflisch R. E.,   Klapper I.  and  Steele G.,
{\it Remarks on singularities, dimension and energy dissipation
for ideal hydrodynamics and MHD}, Comm. Math. Phys., 184 (1997),
443-455.



\bibitem{CCM}  Cannone M.,  Chen Q.-L. and  Miao C.-X.,
{\it A losing estimate for the Ideal MHD equations with application to Blow-up
criterion},  SIAM J. Math. Anal. (2006) DOI.10.1137/060652002.


\bibitem{Ch1}  Chemin J.-Y.,  {\it Perfect incompressible fluids}, Oxford University Press, New York, 1998.

\bibitem{Duv} Duvaut G. and  Lions J.L.  {\it In\'{e}quation  en thermo\'{e}lasticite et magn\'{e}tohydrodynamique},
Arch. Rational Mech. Anal. 46 (1972), 241-279.

\bibitem{Giga}  Giga Y., {\it Solutions for semilinear parabolic equations in $L^p$ and regularity of weak solutions of
the Navier-Stokes system}, J. Diff. Equa., 62(1986), 182-212.


\bibitem{He}  He C. and  Xin Z., {\it On the regularity of weak solutions to the magnetohydrodynamic
                                 equtions}, J. Diff. Equa., 213(2005), 235-254.


\bibitem{KozTa} Kozono  H.  and  Taniuchi Y., {\it Bilinear estimates in BMO and the Navier-Stokes equations}, Math. Z., 235(2000), 173-194.


\bibitem{Koz}  Kozono H.,  Ogawa T. and  Taniuchi Y., {\it The ciritical Sobolev inequalities in Besov spaces and
regularity criterion to some semi-linear evolution equations}, Math. Z., 242(2002), 251-278.

\bibitem{Maj1} Majda A. J.,
{\it Compressible fluid flow and systems of conservation laws in several space
variables},
Applied Mathematical Sciences, 53, Springer-Verlag, New York,
1984.



\bibitem{Nun}  N\'{u}\~{n}ez M., {\it Existence theorems for two-fluid magnetohydrodynamics}, J. Math. Phys., 46(2005), 083101.

\bibitem{Priest}  Priest E.R. and Forbes  T.G.,  {\it Magnetic reconnection: MHD Theory and Applications}.
             Cambridge University Press, Cambridge, UK, 2000.


\bibitem{Ser}  Sermange M. and  Temam R.,  {\it Some mathematical questions related to the MHD equations}, Comm. Pure Appl. Math.,
36(1983), 635-664.



\bibitem{Tri}  Triebel H.,  {\it Theory of Function Spaces}.
             Monograph in mathematics, Vol.78 ,  Birkhauser Verlag,
          Basel, 1983.

\bibitem{Wu0}  Wu J., {\it Bounds and new approaches for  the 3D  MHD equations}, J. Nonlinear Sci.,
                12 (2002), 395-413.

\bibitem{Wu}  Wu J., {\it Regularity results for weak solutions of the 3D  MHD equations}, Discrete Cont. Dyn. S.,
                10 (2004), 543-556.

\bibitem{ZL}  Zhang Z. and  Liu X., {\it On the blow-up criterion of smooth solutions
                to the 3D Ideal MHD equations}, Acta Math. Appl. Sinica, E,
                20 (2004), 695-700.

\bibitem{Zhou}  Zhou Y., {\it Remarks on the regularity for the
                 3D  MHD equations}, Discrete Cont. Dyn. S.,
                12 (2005), 881-886.


\end{thebibliography}
\end{document}